\input amstex
\documentstyle{amsppt}
\NoBlackBoxes
\topmatter
\title Fixed Point Indices and Manifolds with Collars
\endtitle
\author Chen-Farng Benjamin and Daniel Henry Gottlieb
\endauthor
\document
\bigskip

\abstract This paper concerns a formula which relates the Lefschetz 
number $L(f)$ for a
map $f:M\to M'$ to the fixed point index $I(f)$ summed with the fixed point
index of a derived map on part of the boundary of $\partial M$. Here $M$ is a
compact manifold and $M'$ is $M$ with a collar attached.
\endabstract 
\endtopmatter
\medskip
\centerline{\bf 1. INTRODUCTION}

\bigskip\bigskip

This Paper represents the first third of a Ph.~D thesis [{\bf{16}}] 
written by the
first author under the direction of the second author at Purdue
University in 1990. The Thesis is entitled " Fixed Point Indices and
Transfers, and Path Fields" and it 
contains, in addition to the contents of this
manuscript, a formula  analogous to (1.1), which relates to Dold's fixed point 
transfers
and a study of path fields of differential manifolds in
order to relate the formula in this manuscript with an analogous formula
involving indices of vector fields. 

Let $M$ be a compact differentiable
manifold with or without boundary $\partial M$. Assume $V$ is a vector
field on $M$ with only isolated zeros. If $M$ is with boundary $\partial
M$ and $V$ points outward at all boundary points, then the index of the
vector field $V$ equals Euler characteristic of the manifold $M$. This
is the classical Poincar\'e-Hopf Index Theorem. (A 2-dimensional version
of this theorem was proven by Poincar\'e in 1885; in full generality the
theorem was proven by Hopf [{\bf{11}}] in 1926). In particular, the
index is a topological invariant of $M$; it does not depend on the
particular choice of a vector field on $M$. \smallskip Marston Morse
[{\bf{13}}] extended this result to vector fields under more general
boundary conditions, namely, to any vector field without zeros on the
boundary $\partial M$; he discovered the following formula:
$$\text{Ind}(V)+\text{Ind}(\partial\_V)=\chi(M), \tag{1.1} $$ where
$\chi(M)$ denotes the Euler characteristic of $M$ and $\partial\_V$ is
defined as follows: Let $\partial\_M$ be the open subset of the boundary
$\partial M$ containing all the points $m$ for which the vectors $V(m)$
point inward; and let $\partial V$ be the vector field on the boundary
$\partial M$ obtained by first restricting $V$ to the boundary and then
projecting $V|_{\partial M}$ to its component field tangent to the
boundary. Then $\partial\_V=\partial V|_{\partial\_ M}$. Furthermore, in
the same paper, Morse generalized his result to indices of vector fields
with nonisolated zeros. This is the formula (1.1).Now (1.1) was rediscovered by
D.~Gottlieb [{\bf{7}}] and C.~Pugh [{\bf{15}}]. D.~Gottlieb further
found further interesting applications in [{\bf{8}}],[{\bf{9}}] and
[{\bf{10}}]. Throughout this paper we shall call formula (1.1) the Morse 
formula for indices of vector fields.

We consider maps $f:M\to M'$ from a compact topological manifold $M$ to $M'$ 
where $M'$ is obtained by attaching a collar $\partial M\times [0,1]$ to $M$. 
If $f$ has no fixed points on the boundary $\partial M$, 
we prove Theorem (2.2.1) 
which is the fixed point version of the Morse formula:
$$
I(f)+I(r \circ f|_{\partial\_M})=L(r \circ f),
\tag{1.2}
$$
where $I$ denotes the fixed point index, $r$ is a retraction of $M'$ 
onto $M$ which maps the collar $\partial M\times I$ onto the boundary 
$\partial M$, $\partial\_ M$ is an open subset of $\partial M$ containing 
all the points $x\in \partial M$ mapped outside of $M$ under $f$ and $L(r \circ f)$ 
is the Lefschetz number of the composite map $r \circ f$.

In particular, if the map $r \circ f$ is homotopic to the identity map, we have:
$$
I(f)+I(r \circ f|_{\partial\_M})=\chi(M),
\tag{1.3}
$$ 
which is similar to the Morse formula; and the map $r \circ f|\partial\_M$ is analogous to the vector field $\partial\_V$.

Formula(1.3) was independently obtained by A.~Dold (private letter to D.~Gottlieb).

This paper is organized as follows:

In Section  2, \S{1}, we list some properties of fixed point indices; our first main result, Theorem (2.2.1), is proven in Section 2, \S{2}.

\vfill\eject

\medskip
\centerline{\bf 2. FIXED POINT VERSION OF THE MORSE FORMULA}

\bigskip\bigskip

In this section, we use the definition of fixed point index and some well known results on fixed point index given by Dold in [{\bf{3}}] or [{\bf{5}},Chap. 7] to obtain an equation for fixed point indices [Theorem 2.2.1] analogous to the Morse equation for vector field indices described in the introduction.

\medskip

\S{\bf{1. Fixed point index and its properties.}}

\medskip

Let $X$ be an Euclidean neighborhood retract $(ENR)$. Consider maps $f$ from an open subset $V$ of $X$ into $X$ whose fixed point set 
$F(f)=\{x\in V|f(x)=x\}$ is compact. 
Dold [{\bf{3}}] defined the fixed point index $I(f)$ and 
proved the following properties.
\medskip
(2.1.1) LOCALIZATION.
\medskip
Let $f:V\to X$ be a map such that $F(f)$ is compact, then $I(f)=I(f|_W)$ for any open neighborhood $W$ of $F(f)$ in $V$.
\medskip
(2.1.2) ADDITIVITY.
\medskip

Given a map $f:V\to X$ and $V$ is a union of open subsets $V_j$, $j=1,2\dots,n$ such that the fixed point sets $F^j=F(f)\cap V_j$ are mutually disjoint. Then for each $j$, $I(f|_{V_j})$ is defined and 
$$
I(f)=\sum^n_{j=1} I(f|_{V_j}).
$$

(2.1.3) UNITS.
\medskip

Let $f:V\to X$ be a constant map. Then
$$
\align I(f)&=1 \qquad \text{if}\, f(V)=p\in V\quad \text{and}\\
I(f)&=0 
\qquad \text{if}\, f(V)=p\notin V
\endalign
$$

(2.1.4) NORMALIZATION.
\medskip
If $f$ is a map from a compact $E N R \quad X$ to itself, then $I(f)=L(f)$ where $L(f)$ is the Lefschetz number of the map $f$.
\medskip
(2.1.5) MULTIPLICATIVITY
\medskip
Let $f:V\to X$ and $f':V'\to X'$ be maps such that the fixed point sets $F(f)$ and $F(f')$ are compact, then fixed point index of the product $f\times f': V\times V'\to X\times X'$ is defined and
$$
I(f\times f')=I(f)\cdot I(f').
$$

(2.1.6) COMMUTATIVITY.
\medskip
If $f:U\to X'$ and $g:U'\to X$ are maps where $U\subseteq X$ and $U'\subseteq X'$ are open subsets, then the two composites $gf:V=f^{-1}(U')\to X$ and $fg:V'=g^{-1}(U)\to X'$ have homeomorphic fixed point sets. In particular, $I(fg)$ is defined if and only if $I(gf)$ is defined, in that case, 
$$
I(fg)=I(gf).
$$         	

(2.1.7) HOMOTOPY INVARIANCE
\medskip

Let $H:V\times I\to X$ be a homotopy between the maps $f_0$ and $f_1$. Assume the set $F=\{x\in V|H(x,t)=x$ for some $t\}$ is compact; then
$$
I(f_0)=I(f_1).
$$

For our purposes it is useful to reformulate the properties of additivity (2.1.2) and homotopy invariance (2.1.7) in the form of the following propositions. These reformulations are found in R. F Brown's book
[2b], and they form part of an axiom system for the fixed point index. The five axioms are a subset of
Dold's properties. They consist of localization, homotopy invariance , addititvity, normalization and commutivity. We will show that the main formula will follow from these axioms. We will give an alternate proof in the next subsection \S 2.
\medskip
(2.1.8) PROPOSITION.
\medskip
{\it{Assume $X$ is compact and $V$ is an open subset of $X$. Let $f:\overline{V}\to X$ be a map without fixed points on $Bd(V)$. If $\{V_j\}$, $j=1,2,\dots,n$ are mutually disjoint open subsets of $V$ and whose union $\overset n\to{\underset j=1 \to\cup} V_j$ contains all the fixed points of $f$, then}}
$$
I(f|_V)=\sum^n_{j=1}I(f|_{V_j}).
$$

(2.1.9) PROPOSITION.
\medskip
{\it{Assume $X$ is compact and $V$ is an open subset of $X$. Let $\overline H:\overline V\times I\to X$ be a homotopy from $\overline f_0$ and $\overline f_1$ where $\overline f_0$ and $\overline f_1$ are maps from $\overline V$, the closure of $V$ to $X$. If $\overline H(x,t)\neq x$ for all $x\in Bd(V)$ and for all $t$, then 
$$
I(f_0)=I(f_1)\qquad  where\quad f_0=\overline f_0|V\quad and\; f_1=\overline f_1|V.
$$}}

{\it{Proof}}: Since $H=\overline H|_{V\times I}$ is a homotopy from $f_0$ to $f_1$, it suffices to verify that the set $F=\{x\in V|H(x,t)=x$ for some $t\}$ is compact. Let $\{x_j\}$ be a sequence in $F$ converging to $x\in \overline V=V\cup Bd(V)$. There exists a subsequence $\{t_j\}$ of those $t$'s in $I$ such that $\overline H(x_j,t_j)= x_j$. Since $I$ is compact, a  subsequence of $\{t_j\}$ converges to a point $t\in I$. By the continuity of $\overline H$, we have $\overline H(x,t)=x$. On the other hand, we know that $\overline H(x,t)\neq x$ for all $x\in Bd(V)$; thus, $x\in V$ and $H(x,t)=x$. Consequently, $x\in F$. Therefore, $F$ is a closed subset of a compact space, hence $F$ is compact. This proves the proposition.

\bigskip

\S{\bf{2. The main formula.}}
\medskip

Consider a compact topological manifold $M$ with boundary $\partial M$. We attach a collar to $M$ and call the resulting manifold $M':M'=M\underset\partial M\sim\partial M\times\{0\}\to \cup \partial M\times [0,1]$. Let $f:M\to M'$ be a map such that $f(x)\neq x$ for all $x\in \partial M$. Since $M$ is compact, the fixed point set $F(f)$ is a compact set contained in $\overset\circ \to M =M \backslash \partial M$. For such $f:M\to M'$, we define the index of $f$, denoted by $I(f)$, to be the fixed point index of the map $f|_{\overset\circ \to M}$ given in \S{1}.

For specificity we define the retraction $r$:
Let $r:M'\to M$ be the retraction from $M'$ to $M$ given by the formula:
$$
\cases r(m)=m &\text{for $m\in M$}\\
r(b,t)=(b,0)\sim b &\text{for $(b,t)\in\partial M\times [0,1]$.}\endcases
$$
Now we can formulate the main result of the section.

Now, assume $r'$ is any retraction from $M'$ to $M$ such that $r'$ maps the collar $\partial M\times[0,1]$ into the boundary $\partial M$. Then the following Theorem is true:
\medskip
(2.2.1) THEOREM.
$$
I(f)+I(r'f|_{\partial\_ M})=L(r'f).
$$
{\it Furthermore,}
$$\split
L(rf)&=L(r'f)\quad \text{and}\\
I(rf|_{\partial\_ M})&=I(r'f|_{\partial\_ M})
\endsplit
$$
{\it where $r$ is the standard retraction defined above and
where $L(rf)$ denotes the Lefschetz number of $rf:M\to M$ and $\partial \_M=\{x\in\partial M|f(x)\notin M\}$.}
\medskip
{\it Proof}: First, we prove the formula: $I(f)+I(r'f|_{\partial\_ M})=L(r'f)$.
\medskip
Let $V_1=\{x\in M|f(x)\in \overset\circ \to M\}$ and $V_2=\{x\in M|f(x)\in M'\backslash M\}$, then $V_1$ and $V_2$ are disjoint open subsets of the manifold $M$ and $V_1\cup V_2$ contains all the fixed points of the map $r'f$. Indeed, if $x\notin(V_1\cup V_2)$, then $f(x)\in\partial M$ and hence $r'f(x)=f(x)\neq x$. Proposition (2.1.8) implies the equation:
$$
I(r'f)=I(r'f|_{V_1})+I(r'f|_{V_2}).
$$

Since $r'f$ is a self map from $M$ to $M$, so
$$
I(r'f) = L(r'f).
$$
We have:
$$
L(r'f)=I(r'f|_{V_1})+I(r'f|_{V_2}).
\tag{i}
$$
Now, since $r'f|_{V_1}=f|_{V_1}$ and $F(f)\subseteq V_1$, then
$$
I(r'f|_{V_1})=I(f|_{V_1})=I(f).
\tag{ii}
$$
Let us decompose the map $r'f|_{V_2}$:
$$
r'f|_{V_2}:V_2\overset f|_{V_2}\to\longrightarrow \partial M\times [0,1]\overset r'\to\longrightarrow \partial M \overset i\to\longrightarrow M.
$$
The commutativity (2.1.6) implies that
$$
I(r'f|_{V_2})=I(ir'f|_{V_2})=I(r'fi|_{i^{-1}(V_2)})=I(r'f|_{\partial\_ M}).
\tag{iii}
$$
Combining equations (i), (ii) and (iii), we obtain
$$
I(f)+I(r'f|_{\partial\_ M})=L(r'f).
\tag{iv}
$$
This completes the proof of the formula holding for any retraction $r'$. The following two lemmas
will show that the terms in equation (iv) are the same no matter which retraction $r'$ is chosen.

(2.2.2) LEMMA. {\it The retraction $r$ is homotopic to $r'$}
\medskip
{\it Proof}: Consider the homotopy $H_t:M'\to M$, $0\leq t\leq 1$, defined as follows:
$$\split
H_t(m)=m \qquad &\text{for $m\in M$,}\\
H_t(b,s)=r'(b,st) \qquad &\text{for $(b,s)\in \partial M\times [0,1]$.}
\endsplit
$$
Clearly, $H_0=r$ and $H_1=r'$. So, $rf$ and $r'f$ are homotopic.

\medskip

(2.2.3) LEMMA. {\it $L(r'f)=L(rf)$ and $I(rf|_{\partial\_ M})=I(r'f|_{\partial\_ M})$.}
\medskip
{\it Proof}: By LEMMA (2.2.2) $rf$ and $r'f$ are homotopic and, consequently,
$$
L(rf)=L(r'f).
\tag{v}
$$
since the Lefschetz number $L$ is a homotopy invariant.

Equations (iv) and (v) and equation (v) with $r$ replacing $r'$ imply that
$$
I(rf|_{\partial\_ M})=I(r'f|_{\partial\_ M}).
$$
This concludes the proof of 
THEOREM (2.2.1).

\medskip

(2.2.4) COROLLARY. {\it{If $f:M\to M'$ is map such that $f(x)\notin M$ for any $x\in\partial M$, then $I(f)=L(rf)-L(rf|_{\partial M})$.}}
\medskip
(2.2.5) COROLLARY. {\it{If $f:M\to M'$ is without fixed points on the boundary $\partial M$ and $f(\partial M)\subset M$, then $I(f)=L(rf)$.}}
\medskip
(2.2.6) EXAMPLE. {\it{Consider a map $f:D^n\to {\Bbb R}^n$, Here $D^n$ is the unit 
ball and $S^{n-1}$ is the unit boundary sphere, so we can think of $ {\Bbb R}^n $ as $D^n$ with
an open collar attached.

\medskip
\itemitem{\rm (i)} if $f(S^{n-1})\subset D^n$, then $f$ has a fixed point.

\itemitem{\rm (ii)} if $f(S^{n-1})\subset {\Bbb R}^n \backslash D^n$, then Corollary (2.2.4) implies that\newline
 $I(f)=L(rf)-L(rf|_{S^{n-1}})
=1-(1+(-1)^{n-1}$ \rm{deg}$(rf|_{S^{n-1}}))=(-1)^n$\rm{deg}$(rf|_{S^{n-1}})$.}}
\medskip
(2.2.7) COROLLARY. {\it{If $f:M\to M'$ is homotopic to the inclusion map $M\hookrightarrow M'$, then $I(f)+I(rf|_{\partial\_ M})=\chi (M)$ where $\chi (M)$ denotes the Euler characteristic of $M$.}}
\medskip
{\it Proof}: If $f:M\to M'$ is homotopic to the inclusion map $M\hookrightarrow M'$, then the composite map $rf:M\to M$ is homotopic to the identity map. Therefore $L(rf)=L(Id)=\chi (M)$.

(2.2.8) REMARK: {\it Here is a more geometric proof of the Main Theorem (2.2.1)}.
\medskip
{proof}: 
\medskip
Let $DM$ be the double of $M$. That is the union of two copies of $M$ intersecting on their
boundaries. Let $R: DM \rightarrow M$ be the retraction which takes the second copy onto the
first. Now $f \circ R : DM \rightarrow M$. Then the Lefschetz numbers $L(f) = L(f \circ R )$ since 
$R$ is a retraction, which splits the homology of $DM$, so that the traces of the induced map
must be calculated only on the first copy $M$ of $DM$.

Also we consider $M \subset M' \subset DM$. Then $R$ restricted to $M'$ is equal to $r$. Now the
fixed point set of $f \circ R $ consists of the fixed point set of $f$, in the interior of $M$, and the fixed
point set $F(f\circ R) = F (f \circ r)$ contained in $\partial\_ M$. Now the index of $r \circ f$ calculated
on the open set $\partial_ M$ is equal to the index calculated on a small open set $V$ of $M'$ containing $\partial\_  M$ which follows from the next lemma.

(2.2.9) LEMMA.
$$
I(r \circ f|_{\partial\_ M})=I(r \circ f|_V).
$$

{\it{Proof:}} The commutativity (2.1.6) implies that
$$
I(r \circ f|_V)=I(f \circ r|_{r^{-1}(V)}).
$$

It is easy to see that the fixed point set of the map $f \circ r|_{r^{-1}(V)}$ is $\{(b,t)\in\partial M\times (0,1]|f(b)=(b,t)\}$ and the fixed point set of the map $rf|_V$ is $\{b\in\partial M|f(b)=(b,t)$ for some $t\}$.

We now define a homotopy $G_s$, $0\leq s\leq 1$, as the composite of the following maps
$$
\overline{\partial\_ M}\times I \overset r \to\longrightarrow \overline{\partial\_ M}\overset f \to\longrightarrow\partial M\times I\overset H_s \to\longrightarrow \partial M\times I
$$
where the map $H_s$ is defined as follows:
$$
H_s(b,t)=(b,s\overline t +(1-s)t),\qquad \text{where $\overline t$ is a constant, $0<\overline t\leq 1$}.
$$
since the map $H_0$= Identity, we have
$$\split
G_0(x,t)&=H_0(fr(x,t))=fr(x,t),\qquad \text{and}\\
G_1(x,t)&=H_1(fr(x,t))=(rf\times g)(x,t)
\endsplit
$$
where $r \circ f$ is a map from $\overline{\partial\_ M}$  to $\partial M$ and $g:I\to I$, $g(t)=\overline t$, is the constant map. Furthermore, the restriction $G_s|_{Bd(\partial\_ M\times I)}$ has no fixed points for any $0\leq s\leq 1$. To see this, we look at a point $x\in Bd(\partial\_ M)$. We know then $f(x)\in \partial M$ and $rf(x)=f(x)\neq x$, therefore, $G_s(x,t)=H_s(fr(x,t))=H_s(f(x))=H_s(f(x),0)=(f(x),s\overline t)\neq (x,t)$.

Now the properties (2.1.9), (2.1.5) and (2.1.3) imply
$$\split
&I(fr|_{\partial\_ M\times(0,1]})=I(rf|_{\partial\_ M})\cdot I(g)=I(rf|_{\partial\_ M})\cdot 1\\
&I(rf|_V)=I(fr|_{r^{-1}(V)})=I(fr|_{\partial\_ M\times (0,1]}).
\endsplit
$$

The last equality holds because $\partial\_ M\times (0,1]$ contains the fixed point set of $(fr|_{r^{-1}(V)})$.

Thus, $I(rf|_V)=I(rf|_{\partial\_ M})$.
\medskip
{\it{Proof of Theorem (2.2.1):}}
\medskip
Consider the composite $M\overset f\to \longrightarrow M'\overset r\to \longrightarrow M$. Let $V$ be the open set as in Lemma (2.2.3), then $V$ and $\overset\circ \to M$ are two open subsets of $M$ such that $V\cup\overset\circ \to M=M$. Clearly, $F(rf)\cap\overset\circ \to M$ and $F(rf)\cap V$ are disjoint. Using additivity (2.1.2) and the normalization (2.1.4) of the fixed point indices, we have:
$$
I(rf|_{\overset\circ \to M})+I(rf|_V)=I(rf)=L(rf).
$$
Lemma(2.2.2) and Lemma (2.2.3) then imply the equation:
$$
I(f)+I(rf|_{\partial\_ M})=L(rf).
$$

\vfill\eject

\eject
\enddocument
\bye